\documentclass[letterpaper, preprint, paper,11pt]{AAS}	

\usepackage{bm}
\usepackage{amsmath}
\usepackage{subfigure}
\usepackage[colorlinks=true, pdfstartview=FitV, linkcolor=black, citecolor= black, urlcolor= black]{hyperref}
\usepackage{overcite}
\usepackage{footnpag}			      	

\PaperNumber{23-181}

\begin{document}

\title{CLOSED-FORM OPTIMAL PROPULSIVE-DIFFERENTIAL DRAG CONTROL FOR LARGE RECONFIGURATIONS OF SPACECRAFT SWARMS}

\author{Matthew Hunter\thanks{Ph.D. Candidate, Stanford University, Department of Aeronautics \& Astronautics, Space Rendezvous Laboratory, Durand Building, 496 Lomita Mall Stanford, CA 94305}
\ and Simone D'Amico\thanks{Assistant Professor, Stanford University, Department of Aeronautics \& Astronautics, Space Rendezvous Laboratory, Durand Building, 496 Lomita Mall Stanford, CA 94305}
}

\maketitle{}

\begin{abstract}
Distributed Space Systems offer unique capabilities but have significant hardware limitations on computational processing and fuel when employing small satellites. These constraints drive interest in deterministic maneuver planning that incorporates orbital perturbations to improve final state accuracy and decrease reconfiguration delta-v cost. This paper proposes a novel closed-form, provably optimal framework to generate complete roto-translational solutions for propulsive-differential drag control of multiple satellites. The approach leverages relative orbit elements and reachable set theory to characterize all control with a single methodology that can be readily generalized to generic impulsive-continuous control. The algorithm is robust to large reconfigurations and extended control windows where differential drag is most effective.
\end{abstract}

\section{INTRODUCTION}

Distributed Space Systems (DSS) have drawn increasing interest from the space community as their unique capabilities, such as highly variable instrument focal length, robustness through redundancy, and in-flight servicing, have been successfully demonstrated on PRISMA, TanDEM-X, Orbital Express, and the Mission Extension Vehicle \cite{damico_phd,DSS_survey}. Using smaller, often modular satellites introduces additional complex constraints, chief among them limits in processing power and fuel. The advantages and constraints of DSS are evident on the Space Weather Atmospheric Reconfigurable Multiscale Experiment (SWARM-EX) mission, which will use a novel swarm of three 3U CubeSats in Low Earth Orbit (LEO) to take distributed ionospheric and thermospheric space weather measurements, requiring multiple changes of nearly 1000km in inter-spacecraft separation \cite{swarm-ex}. Such a mission is not feasible under the limitations of CubeSat propellant, motivating interest in both optimal propulsive maneuvering and the use of environmental forces whose differential effect between spacecraft can be controlled. 

Propulsive maneuvers are often modelled as "impulsive", or instantaneous additions of delta-v, due to their relatively high magnitude and short duration with respect to orbital period. Extensive literature exists for impulsive control, with approaches that vary widely in applicability and computational efficiency. The nonlinear dynamics of cartesian relative spacecraft motion often necessitate numerical optimization techniques such as sequential convex programming, Lie bracket theory, sliding mode control, and differential dynamic programming \cite{guffanti_SCP,sipowa_Liebracket,cho_robustImp,sasaki_DDP}. Without proper constraints, numerical methods often converge to kinematically sub-optimal continuous control profiles. Integration Constant (IC) states, such as Relative Orbital Elements (ROE), can simplify the dynamics to linear time variant (LTV) models with minimal loss of accuracy, enabling closed-form maneuver planning on missions like PRISMA and TanDEM-X \cite{damico_phd}. In general, impulsive optimal control problems with norm-like cost functions and LTV dynamics can be solved efficiently through reachable set theory and, when combined with an ROE parameterization, can produce provably optimal closed-form solutions for orbits of arbitrary eccentricity under J2-perturbed dynamics \cite{koenig_impulse,chernick}.

The main perturbative force in LEO behind J2 is atmospheric drag, the key motivation of this work, and modulating differences in DSS spacecraft attitude produces differential drag. Differential drag cannot usually achieve full state controllability due to strict limits on its magnitude and directionality and poor on-board characterization of atmospheric density. A desire to guarantee final state accuracy has resulted in many numerical techniques to incrementally approach a desired formation without thrust through Linear Quadratic Regulator (LQR), Model Predictive Control (MPC), Proportional Integral Derivative (PID), or nonlinear switching lane formulations \cite{bev_drag_roto-trans,bev_drag_leastsquares,harris_LQRdrag,ben-larbi_drag_ROE,koenig_drag}.  Focusing instead on the roto-translational problem, simultaneous attitude stabilization and translational state changes have been accomplished with dual-quaternion parameterizatons and on spacecraft with a strictly specialized drag construction, all with provable Lyapunov stability \cite{stanfield_dualq_lqr,yang_dualq}. Overall, these control efforts have primarily limited themselves to considering propulsive or differential drag control alone, often lacking the computational efficiency to be implemented on-board.

In contrast, this paper proposes a provably fuel-optimal, closed-form architecture for LEO prop-ulsive-differential drag control, also known as hybrid control. The hybrid approach can decrease the delta-v cost of swarm reconfigurations, as drag maneuvers have little or no propellant cost and higher-thrust propulsive maneuvers provide the control authority drag lacks. This novel framework is realized through reachable set theory, demonstrating the ability of a single methodology to produce closed-form provably optimal control solutions for multiple methods of control, and results in two major contributions to the state of the art. First, the full hybrid architecture represents the first ever complete framework for provably optimal hybrid DSS maneuver planning, in an approach that can be easily extrapolated to generic impulsive-continuous LTV control systems. Second, the reachable set theory analysis of hybrid control is leveraged to characterize optimality, increase final state accuracy, and improve computational efficiency for extended control windows and large reconfigurations in near-circular orbit. The hybrid architecture is presented in three major sections. First, reachable set theory identifies optimal maneuver times and directions for both types of actuation for the overall hybrid problem. Second, algorithms are designed to produce roto-translational maneuver plans for reconfigurations of arbitrary size and duration, with the computational efficiency to be used on-board CubeSats. Finally, optimality, computational efficiency, and final state accuracy are verified against numerical and closed-form approaches in a high fidelity full-force orbital simulation that propagates the Gauss Variational Equations (GVE) for all relevant LEO perturbations.

\section{PROBLEM DEFINITION}
\subsection{State Representation}

This paper considers the relative motion between a chief (a non-maneuvering spacecraft or reference orbit) and a deputy spacecraft, indicated by the subscripts $c$ and $d$ respectively. The relative motion and control actions of the deputy with respect to the chief are commonly represented in cartesian frames, primarily Radial-Tangential-Normal (RTN). RTN is defined from the spacecraft's center of mass and points radially outward from the Earth's center of mass, normally along the angular momentum vector, and tangentially in the direction of motion to complete the right-handed triad. Spacecraft relative motion can be equivalently expressed in ROE, an advantageous IC state, to approximate nonlinear cartesian dynamics through linear models \cite{damico_phd}. This work uses quasi-nonsingular ROE, valid for inclined orbits of arbitrary eccentricity, defined as
\begin{equation} \label{eq: ROE}
    \delta\boldsymbol{\alpha} = 
    \begin{bmatrix}
    \delta a \\ \delta \lambda \\ \delta e_x \\ \delta e_y \\ \delta i_x \\ \delta i_y
    \end{bmatrix} = 
    \begin{bmatrix}
    \Delta a/a_c \\
    \Delta u + \Delta \Omega \cos i_c \\
    e_d\cos\omega_d - e_c\cos\omega_c \\
    e_d\sin\omega_d - e_c\sin\omega_c \\
    \Delta i \\
    \Delta \Omega \sin i_c
    \end{bmatrix}
\end{equation} 
\noindent where $\Delta$ is a difference in the associated quantity between the chief (subscript $c$) and deputy (subscript $d$), $\delta a$ is the relative semi-major axis, $\delta \lambda$ is the relative mean longitude, $\delta e_x$ and $\delta e_y$ are the components of the relative eccentricity vector $\delta \boldsymbol{e}$, $\delta i_x$ and $\delta i_y$ are the components of the relative inclination vector $\delta \boldsymbol{i}$, $u = M + \omega$ is the mean argument of latitude, and $a$, $e$, $i$, $\Omega$, $\omega$, and $M$ are the Keplerian orbital elements. Mean ROE, derived from the mean Keplerian orbital elements, are constant and equivalent to the IC of the Hill-Clohessy-Wiltshire (HCW) and Yamanaka-Ankersen (YA) differential equations of spacecraft relative motion at small separations in unperturbed Keplerian orbit and vary slowly in the presence of perturbations. A State Transition Matrix (STM) can incorporate the primary perturbations to Keplerian orbit and their corresponding effects on ROE in closed-form and propagate them over a desired time interval. This work uses an STM that includes corrections for J2 for mean ROE in near-circular orbit, given as
\begin{equation} \label{eq: Gamma matrix}
    \boldsymbol{\Phi}^{J_2}(t) = \boldsymbol{\Phi}^{J_2} (\boldsymbol{\alpha}_c(t),t_f-t) = 
    \begin{bmatrix}
    1 & 0 & 0 & 0 & 0 & 0 \\
    \Phi_{21}^{J_2} & 1 & 0 & 0 & \Phi_{25}^{J_2} & 0 \\
    0 & 0 & \cos(\Dot{\omega}_c(t_f-t)) & -\sin(\Dot{\omega}_c(t_f-t)) & 0 & 0 \\
    0 & 0 & \sin(\Dot{\omega}_c(t_f-t)) & \cos(\Dot{\omega}_c(t_f-t)) & 0 & 0 \\
    0 & 0 & 0 & 0 & 1 & 0 \\
    \Phi_{61}^{J_2} & 0 & 0 & 0 & \Phi_{65}^{J_2} & 1 \\
    \end{bmatrix} 
\end{equation}
\noindent with the following simplifying substitutions
\begin{align}
    \eta_c = \sqrt{1-e_{c}^{2}} \text{ , } \kappa_c = \frac{3}{4}\frac{J_2 R^{2}_{E} \sqrt{\mu}}{a_{c}^{7/2} \eta_{c}^{4}}& \text{ , } \Dot{\omega}_c = \kappa_c(5\cos^2(i_c)-1) \\
    \Phi_{21}^{J_2} = -(\frac{3}{2}n_c+\frac{7}{2}\kappa_c(1+\eta_c)&(3\cos^2(i_c)-1))(t_f-t) \nonumber \\
    \Phi_{25}^{J_2} = -\kappa_c(4+3\eta_c)\sin(2i_c)(t_f-t) \text{ , } \Phi_{61}^{J_2} = \frac{7}{2}\kappa_c&\sin(2i_c)(t_f-t) \text{ , } \Phi_{65}^{J_2} = 2\kappa_c\sin^2(i_c)(t_f-t) \nonumber
\end{align}
\noindent where $\boldsymbol{\alpha}_c(t)$ are the mean Keplerian orbital elements at time $t$, $t_f$ is the end of the control window, $n$ is the mean motion, and $\Phi^{J_2}_{i,j} (t,t_f)$ is the element in the $i$th row of the $j$th column \cite{koenig_STM}. The near-circular assumption naturally follows from the focus on LEO spacecraft and their minimally eccentric orbits.

The control input matrix $\boldsymbol{B}$ transforms RTN maneuvers to changes in osculating or approximately mean ROE in near-circular orbit, defined as
\begin{equation} \label{eq: B matrix}
    \boldsymbol{B}(t) = 
    \frac{1}{a_c n_c}
    \begin{bmatrix}
        0 & -2 & \sin u_c & -\cos u_c & 0 & 0 \\
        2 & 0 & 2\cos u_c & 2\sin u_c & 0 & 0 \\
        0 & 0 & 0 & 0 & \cos u_c & \sin u_c \\
    \end{bmatrix}^{\text{T}}
\end{equation}
\noindent where $\boldsymbol{B}(t) = \boldsymbol{B}(\boldsymbol{\alpha}_c(t))$ and $n$ is the mean motion \cite{damico_phd}. It is important to note that $\boldsymbol{B}$ decouples in-plane (IP) control (radial and tangential maneuvers with $\delta a$, $\delta \lambda$, and $\delta \boldsymbol{e}$ states) and out-of-plane (OOP) control (normal maneuvers with $\delta \boldsymbol{i}$ states) as a product of linearization. The dynamics model of relative spacecraft motion can be posed as
\begin{equation} \label{eq: dynamics}
    \Delta \delta \boldsymbol{\alpha} =  \delta \boldsymbol{\alpha}_f - \boldsymbol{\Phi}^{J_2} (t) \delta \boldsymbol{\alpha}_0 = \int_{t_0}^{t_f} \boldsymbol{\Phi}^{J_2} (t) \boldsymbol{B} (t) \boldsymbol{u}(t) dt
\end{equation}
\noindent where $\Delta \delta \boldsymbol{\alpha}$ is the pseudostate of the reconfiguration, defined by initial state $\delta \boldsymbol{\alpha}_0$ at $t_0$ and final desired state $\delta \boldsymbol{\alpha}_f$ at $t_f$, and the RTN actions $\boldsymbol{u}(t)$ occur in the reconfiguration control window $[t_0,t_f]$. It is worth mentioning here that $\boldsymbol{\Phi}^{J_2}$, unlike $\boldsymbol{B}$, produces IP/OOP state coupling that scales linearly with propagation time $t_f-t$, and this interaction will be discussed in later sections. 

\subsection{Differential Drag}

The force of atmospheric drag on a single spacecraft is modelled in RTN as
\begin{equation} \label{eq: atmospheric drag}
    \boldsymbol{p}^{{\text{drag}}} = 
    \begin{bmatrix}
        0 & \frac{1}{2}\rho v^2 \frac{C_D A}{m} & 0
    \end{bmatrix}^{\text{T}} = 
    \begin{bmatrix}
        0 & \frac{1}{2}\rho n^2 a^2 B & 0
    \end{bmatrix}^{\text{T}}
\end{equation}
\noindent where $\rho$ is the atmospheric density, $v = na$ approximates spacecraft tangential velocity in near-circular orbit, $C_D$ is the spacecraft's coefficient of drag, $A$ is the cross-sectional area, $m$ is the spacecraft's mass, and $B = \frac{C_D A}{m}$ defines the spacecraft's ballistic coefficient. The magnitude of atmospheric drag is inversely correlated to spacecraft altitude, such that differential drag is most effective and a viable actuator in LEO. This motivates this work's overall focus on LEO spacecraft and the corresponding near-circular assumption. Differential drag is modelled as a control force in RTN as
\begin{align} \label{eq: drag force}
    \boldsymbol{u}^{\text{drag}}(t) = 
    \boldsymbol{p}^{\text{drag}}_c(t) - \boldsymbol{p}^{\text{drag}}_d(t) = 
    \begin{bmatrix}
    0 & \frac{1}{2}n_c^2a_c^2\Delta B^{\rho}(t) & 0
    \end{bmatrix}^{\text{T}}
\end{align}
\noindent where $\Delta B^{\rho}(t) = \rho_c(t) B_c(t) - \rho_d(t) B_d(t)$ defines the augmented differential ballistic coefficient. $\Delta B^{\rho}$ is primarily bounded by the chief and deputy's minimum and maximum cross-sectional area. Given the nature of formation flight, the chief and deputy are assumed to be in close orbits and share the chief's $n$, $a$, and RTN frame, since RTN orientation varies by less than $1^\circ$ for along-track separations up to 100km. This simplifies the drag force to a tangential-only maneuver. It is noted here that errors introduced by this model can be compensated by MPC, and this approach will be elaborated on briefly in the conclusions.

\subsection{Optimal Hybrid Control Problem}

The most compact description of the reconfiguration control problem addressed by this paper is given by
\begin{align} \label{eq: hybrid control problem}
    \text{minimize  } &\sum_{i=1}^{k} \left\| \boldsymbol{C} \boldsymbol{u}_i \right\|_{2} \\
    \text{subject to  }& \Delta \delta \boldsymbol{\alpha} = \sum_{i=1}^{k} \boldsymbol{\Gamma}(t_{i}) \boldsymbol{u}_{i} \nonumber
\end{align}
\noindent for $t_i \in [t_0,t_f]$. $\Delta \delta \boldsymbol{\alpha}$ is derived from the reconfiguration dynamic constraint in Eq. \ref{eq: dynamics} and is invariant with respect to the control problem. Control inputs $\boldsymbol{u}_i = [\boldsymbol{u}_{i}^{\text{propulsive}};\boldsymbol{u}_{i}^{\text{drag}}]$ concatenate both forms of control in the RTN frame into a single vector. Assuming the spacecraft has a single propulsive thruster and uses no propulsion to control attitude, the cost function sums the L2-norm of each propulsive action, isolated by selector matrix $\boldsymbol{C} = [\boldsymbol{I}_3,\boldsymbol{0}_3]$ where $\boldsymbol{I}_3$ and $\boldsymbol{0}_3$ are the 3rd-order identity and zero matrices. Deviations from the drag maneuver plan for propulsive maneuvering are short in duration with respect to the control windows relevant for hybrid control. Therefore, the roto-translational constraints between propulsion and drag maneuvers are neglected, such that both maneuver types can be conducted independently and simultaneously. The auxiliary matrix $\boldsymbol{\Gamma}(t_i) = [\boldsymbol{\Phi}^{J_2}(t_i)\boldsymbol{B}(t_i),\int_{t_{i}}^{t_{i+1}}\boldsymbol{\Phi}^{J_2}(t)\boldsymbol{B}(t)dt]$ translates propulsive maneuvers, approximated as impulsive maneuvers of unbounded magnitude, at time $t_i$ and analytically integrates differential drag maneuvers held constant over time step $[t_{i},t_{i+1}]$ to pseudostate space. The control window $[t_0,t_f]$ is discretized into $k$ time steps, such that $\sum_{i=1}^{k} (t_{i+1} - t_{i}) = t_f - t_0$.

\section{REACHABLE SET THEORY}

This work uses reachable set theory to identify minimum reconfiguration costs and optimal maneuver times geometrically as a analytical function of spacecraft state. The following summarizes the approach for a single control source. Let $U(c,t)$ be the set of control actions $\boldsymbol{u}$ in the RTN frame with a norm-like cost no greater than $c$ at time $t$. Let $S(c,t)$ be the set of pseudostates $\Delta \delta \boldsymbol{\alpha}$ that can be reached by a single action $\boldsymbol{u} \in U(c,t)$ under the dynamic constraint in Eq. \ref{eq: dynamics}, given as 
\begin{equation} \label{eq: RST single maneuver}
    S(c,t) = \{\Delta \delta \boldsymbol{\alpha}: \Delta \delta \boldsymbol{\alpha} = \boldsymbol{\Phi}^{J_2}(t) \boldsymbol{B}(t) \boldsymbol{u}, \boldsymbol{u} \in U(c,t)\}
\end{equation}
Let $S(c,T)$ be the set of $\Delta \delta \boldsymbol{\alpha}$ that can be reached by a single control action of cost no greater than $c$ within the control window. Finally, let $S^*(c,T)$ be the set of $\Delta \delta \boldsymbol{\alpha}$ that can be reached by a set of $k \geq 1$ control inputs of combined cost no greater than $c$ taken within the control window. Given that the cost of a control action scales linearly with its magnitude and $\boldsymbol{\Phi}^{J_2}(t) \boldsymbol{B}(t) \boldsymbol{u}_j \in S(c,T)$, $S^*(c,T)$ can also be formulated as a linear combination of the $\Delta \delta \boldsymbol{\alpha}_j$ in $S(c,T)$ normalized by the combined cost of the control actions, given as
\begin{equation} \label{eq: RST convex hull}
    S^*(c,T) = \{\Delta \delta \boldsymbol{\alpha}: \Delta \delta \boldsymbol{\alpha} = \sum^{k}_{j=1} \hat{c}_j \Delta \delta \boldsymbol{\alpha}_j, \Delta \delta \boldsymbol{\alpha}_j \in S(c,T), \hat{c}_j \geq 0, \sum^{k}_{j=1} \hat{c}_j = 1 \}
\end{equation}
\noindent which demonstrates that $S^*(c,T)$ is the convex hull of $S(c,T)$. The cost that scales $S^*(c,T)$ such that its boundary contains $\Delta \delta \boldsymbol{\alpha}$ is the minimum cost $\delta v_{\text{min}}$ for the reconfiguration. 

This formulation enables several important observations about optimal control \cite{chernick}. First, a 2$n$-dimensional (2$n$D) reconfiguration can be broken down into $n$ 2D planes, such that the 2$n$D reconfiguration $\delta v_{\text{min}}$ is equal to or greater than the highest $\delta v_{\text{min}}$ within an isolated 2D plane, neglecting the other dimensions. This establishes the concept of a dominant plane, the 2D plane across all $n$ 2D planes that contains the highest $\delta v_{\text{min}}$, and a dominance case, the specific dimension that drives the $\delta v_{\text{min}}$ within the dominant plane. Second, consider the nested reachable set $S_n^*(c,T_{\text{opt}})$, which is the convex hull of $S_n(c,T_{\text{opt}})$, a nested set of $S(c,T)$ restricted to maneuvers at optimal maneuver times $T_{\text{opt}}$. The reconfiguration $\delta v_{\text{min}}$ equals that of the reconfiguration's dominance case only when $\Delta \delta \boldsymbol{\alpha}$ is reachable at the dominance case $\delta v_{\text{min}}$, meaning $\Delta \delta \boldsymbol{\alpha}$ lies on the boundary of $S_n^*(\delta v_{\text{min}},T_{\text{opt}})$ in the dominant plane and within $S_n^*(\delta v_{\text{min}},T_{\text{opt}})$ in the non-dominant planes. Koenig and D'Amico used this formulation to show that the $\delta v_{\text{min}}$ and $T_{\text{opt}}$ of a dominance case are found as the solution to the optimization problem
\begin{equation} \label{eq: delta-v min}
    \text{maximize: } \frac{\boldsymbol{\eta}^{\text{T}} \Delta \delta \boldsymbol{\alpha}}{\max_{t \in T} \left( \max_{\boldsymbol{u} \in U(1,t)} \left( \boldsymbol{\eta}^{\text{T}} \boldsymbol{\Phi}^{J_2}(t) \boldsymbol{B}(t) \boldsymbol{u} \right) \right) }
\end{equation}
\noindent where $\boldsymbol{\eta}$ is the outward normal of the supporting hyperplane for the contour of the reachable set that intersects $\Delta \delta \boldsymbol{\alpha}$ at minimum cost \cite{koenig_impulse}. While this problem often requires numerical methods for the 6D ROE reconfiguration, the simplicity of the reachable sets in 2D ROE space allows the $\boldsymbol{\eta}$ for each dominance case to be characterized geometrically. The dominance case optimal maneuver times, $T_{\text{opt}}$, correspond to points on $S^*$, defined by $\boldsymbol{\eta}$, that share a boundary with $S$. Furthermore, the optimal maneuver direction in pseudostate space matches the direction of $\boldsymbol{\eta}$. This geometric approach for optimal control will be illustrated in the following sections.

\section{REACHABLE SET ANALYSIS}

The reachable sets for hybrid control are found by taking a single maneuver of arbitrary cost and numerically sampling the $\Delta \delta \boldsymbol{\alpha}$ effect over all possible maneuver times, directions, and durations in the control window to find the set $S$ from Eq. \ref{eq: RST single maneuver}. The convex hull of these samples is set $S^*$ from Eq. \ref{eq: RST convex hull} and the reachable set. The ROE are arranged into 3 2D planes: the $\Delta \delta \boldsymbol{a}$ plane ($\Delta \delta a$, $\Delta \delta \lambda$), the $\Delta \delta \boldsymbol{e}$ plane, and the $\Delta \delta \boldsymbol{i}$ plane. Because drag maneuvers have no cost, propulsive control can be analyzed in isolation to determine dominance cases and optimal propulsive maneuver locations. For the periodic dynamics in the $\Delta \delta \boldsymbol{e}$ and $\Delta \delta \boldsymbol{i}$ planes, the phase of $\eta$ most easily relates to a specific $u_c$, which can be converted to a series of maneuver times.

\subsection{Propulsive Reachable Sets}

\begin{figure} [h]
    \centering
    \includegraphics[scale=0.5]{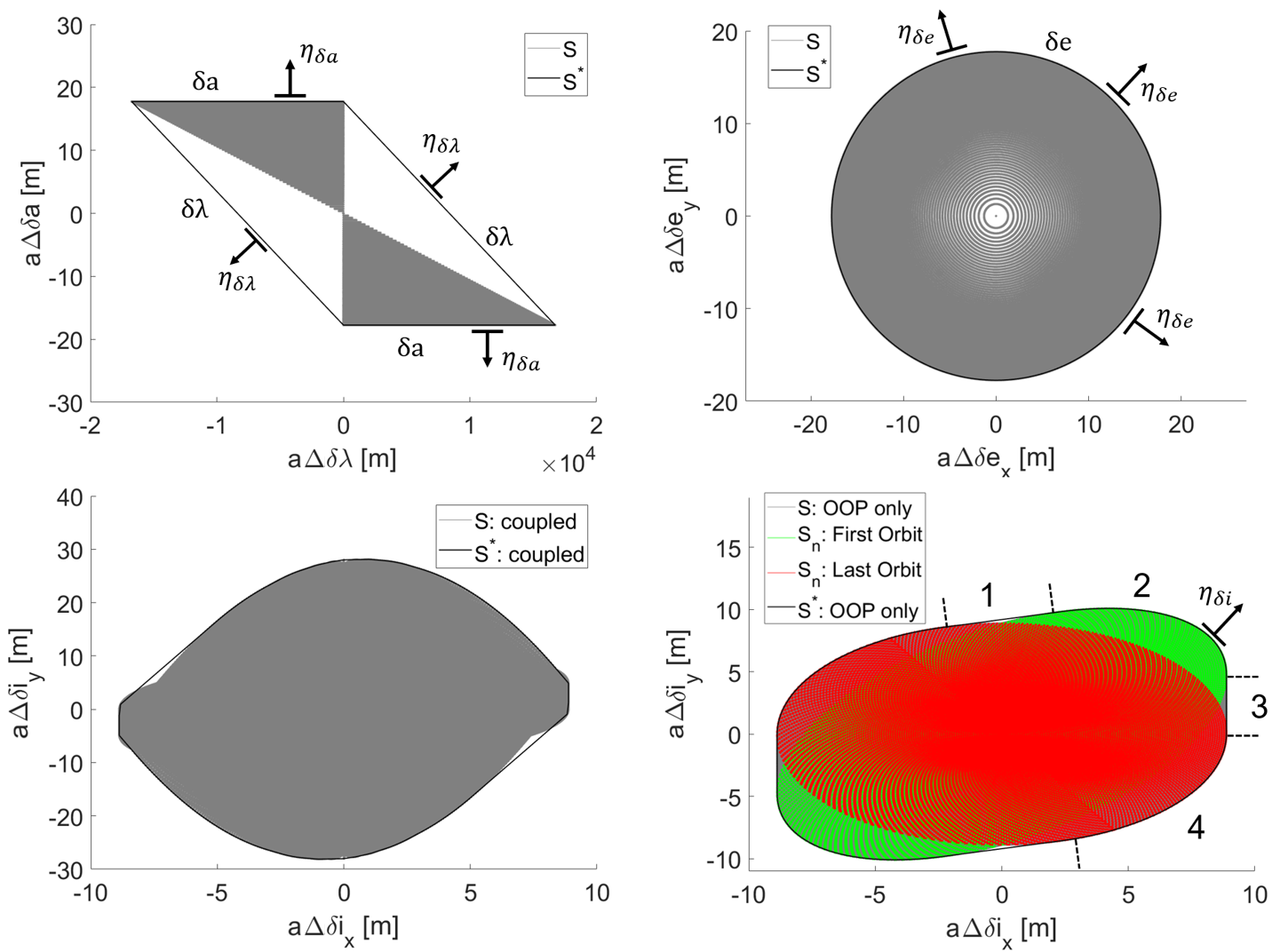}
    \caption{Propulsive reachable sets over a 100 orbit control window for coupled maneuvering in the $\Delta \delta \boldsymbol{a}$ (top left), $\Delta \delta \boldsymbol{e}$ (top right), and $\Delta \delta \boldsymbol{i}$ (bottom left) planes and for decoupled OOP maneuvers in the $\Delta \delta \boldsymbol{i}$ plane (bottom right). Each distinct contour of $S^*$, affine in the $\Delta \delta \boldsymbol{a}$ plane and circular in the $\Delta \delta \boldsymbol{e}$ plane, is labelled with its corresponding dominance case and viable $\eta$.}
    \label{fig: prop reachable sets}
\end{figure}

Figure \ref{fig: prop reachable sets} displays the reachable sets for propulsive control in near-circular orbit, with indicated dominance cases, over an extended 100 orbit control window for visual clarity of the warping due to J2. The IP/OOP coupling from $\boldsymbol{\Phi}^{J_2}$ only warps the reachable set in the $\Delta \delta \boldsymbol{i}$ plane, in the $\Delta \delta i_y$ direction (see Figure \ref{fig: prop reachable sets}, bottom row). The lack of coupling effects in the $\Delta \delta \boldsymbol{a}$ and $\Delta \delta \boldsymbol{e}$ planes implies that $\Delta \delta \boldsymbol{\alpha}$ rarely lies within $S_n^*(\delta v_{\text{min}},T_{\text{opt}})$ in those planes when $\Delta \delta \boldsymbol{i}$ is the dominant plane. Therefore, the 6D control problem is split, according to the partial decoupling in Eq. \ref{eq: B matrix}, into separated 4D IP and 2D OOP sub-problems.

First, the IP dominance cases and related $\eta$ are identified \cite{chernick}. Two distinct affine contours define the $\delta a$ and $\delta \lambda$ dominance cases in the $\Delta \delta \boldsymbol{a}$ plane. $\delta a$ optimal maneuvers can occur anytime in the control window in the direction of $\Delta \delta a$ (see $S^*$ contours matching $\eta_{\delta a}$ in Figure \ref{fig: prop reachable sets}, top left). In contrast, $\delta \lambda$ optimal maneuvers only occur at the beginning and end of the control window in the direction of $\Delta \delta \lambda$ (see $S^*$ contours matching $\eta_{\delta \lambda}$ in Figure \ref{fig: prop reachable sets}, top left). One consistent circular contour exists in the $\Delta \delta \boldsymbol{e}$ plane, the $\delta e$ dominance case (see $S^*$ contours matching $\eta_{\delta e}$ in Figure \ref{fig: prop reachable sets}, top right). In this case, $\delta e$ optimal maneuver locations are restricted by the angle of $\eta$ or phase of $\Delta \delta \boldsymbol{e}$, which occurs at times given by the set
\begin{equation} \label{eq: T opt de prop}
    T_{\text{opt},\delta e} = ( \phi_{\Delta \delta \boldsymbol{e}} + m\pi - \Dot{\omega}_c t_f - u_{c,0} + k_i\pi) / \Dot{M}_c
\end{equation}
\noindent where $\Dot{M}_c = \kappa_c \eta_c (3\cos^2(i_c)-1)$, $u_{c,0}$ is the quantity at $t_0$, $\phi_{\Delta \delta e}=\tan^{-1}(\Delta \delta e_y / \Delta \delta e_x)$, $m$ is an integer such that $\phi_{\Delta \delta e} + m\pi > \Dot{\omega}_c t_f + u_{c,0}$, and $k_i$ is any integer such that $T_{\text{opt},\delta e} \in [t_0,t_f]$.

For OOP control, the entire plane is dictated by the $\Delta \delta \boldsymbol{i}$ vector, but the reachable set over an extended control window differs from that found by Chernick and D'Amico \cite{chernick}. The $\delta i$ dominance case does not have a consistent contour, with the direction of $\eta$ depending on the phase of $\Delta \delta \boldsymbol{i}$. To differentiate these regions, notice that the nonlinear contours of $S^*$ correspond to a subset of $S$ restricted to the first or last orbit of the control window (see Figure \ref{fig: prop reachable sets}, bottom right). Therefore, $\delta i$ dominant regions 2 and 4 are restricted to a single optimal location in the first or last orbit at the phase $\phi_{\Delta \delta \boldsymbol{i}}$, while the affine $\delta i$ dominant regions 1 and 3 have an optimal location in both the first and last orbit at different $\phi_{\Delta \delta \boldsymbol{i}}$. The regions of the $\delta i$ dominance case are summarized in Table \ref{tab: di regions} and correspond to maneuver times given as
\begin{equation} \label{eq: T opt di}
    T_{\text{opt},\delta i} = ( \phi_{\Delta \delta i} + m\pi - \omega_{c,0} + k_i\pi) / (\Dot{M}_c + \Dot{\omega}_c)
\end{equation}

\begin{table} [h]
\caption{$\delta i$ Dominance Case Regions.}
\label{tab: di regions}
\begin{center}
\begin{tabular}{c|c|c}
\hline
Region (Orbit) & Bounds & Optimal $\phi_{\Delta \delta \boldsymbol{i}}$\\
\hline
\hline
1 (First \& Last) & $\tan(u_{c,\delta i}) \leq \frac{\Delta \delta i_y}{\Delta \delta i_x}$ & $u_{c,\delta i} \text{ \& } u_{c,\delta i} + \pi $\\
\hline
2 (First) & $\Phi_{65}^{J_2} \leq \frac{\Delta \delta i_y}{\Delta \delta i_x} \leq \tan(u_{c,\delta i})$ & $\tan^{-1} \left( \frac{\Delta \delta i_y}{\Delta \delta i_x}  - \Phi_{65}^{J_2} \right)$\\
\hline
3 (First \& Last) & $0 \leq \frac{\Delta \delta i_y}{\Delta \delta i_x} \leq \Phi_{65}^{J_2}$ & $\Phi_{65}^{J_2} \text{ \& } 0$\\
\hline
4 (Last) & $-\tan(u_{c,\delta i}) \leq \frac{\Delta \delta i_y}{\Delta \delta i_x} \leq 0$ & $\tan^{-1} \left( \frac{\Delta \delta i_y}{\Delta \delta i_x} \right)$ \\
\hline
\multicolumn{3}{c}{*$u_{c,\delta i} = \cot^{-1} \left( \pm 0.5 \boldsymbol{\Phi}^{J_2}_{65}(t=t_0) \right)$ for notational simplicity}
\end{tabular}
\end{center}
\end{table}

\subsection{Differential Drag Reachable Sets}

\begin{figure} [h]
    \centering
    \includegraphics[scale=0.4]{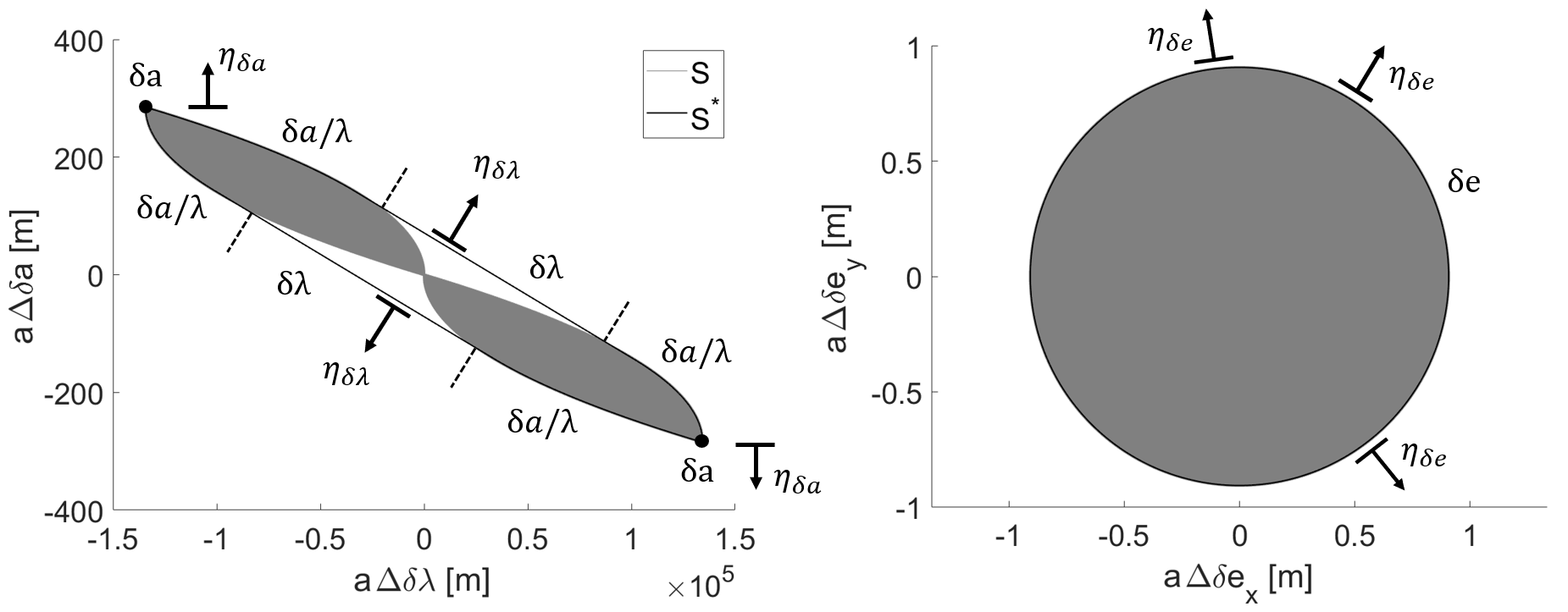}
    \caption{Differential drag reachable sets over a 100 orbit control window under constant atmospheric density in the $\Delta \delta \boldsymbol{a}$ (left) and $\Delta \delta \boldsymbol{e}$ (right) planes. The contours of $S^*$ that match the $\eta$ of the propulsive dominance cases are labelled with those dominance cases and viable $\eta$, along with the regions where the $\delta a$ and $\delta \lambda$ dominance cases have equal $\delta v_{\text{min}}$.}
    \label{fig: drag reachable sets}
\end{figure}

Figure \ref{fig: drag reachable sets} displays the reachable sets for differential drag control in near-circular orbit, with indicated dominance cases. Differential drag is a tangential-only maneuver in Eq. \ref{eq: drag force}, so drag reachable set analysis can be restricted to the IP control sub-problem. The sets assume constant atmospheric density to maintain contour consistency. Equal positive and negative limits on $\Delta B^\rho$ are also assumed for visual clarity, although this assumption is not required for the analysis. While differential drag causes the overall hybrid reachable sets from Eq. \ref{eq: hybrid control problem} to differ from Figure \ref{fig: prop reachable sets}, it introduces no additional cost and does not change the dominance cases or related $\eta$ found for propulsive-only control. These same $\eta$ can be applied to the drag reachable sets to find optimal drag profiles, or the modulation of $\Delta B^\rho$ over time, for each dominance case. The regions in Figure \ref{fig: drag reachable sets} that do not match a propulsive $\eta$ correspond to multiple dominance cases having equal $\delta v_{\text{min}}$ (see $S^*$ contours labelled $\delta a/\lambda$ in Figure \ref{fig: drag reachable sets}, left), a problem that will be addressed in later sections. 

\begin{figure} [h]
    \centering
    \includegraphics[scale=0.42]{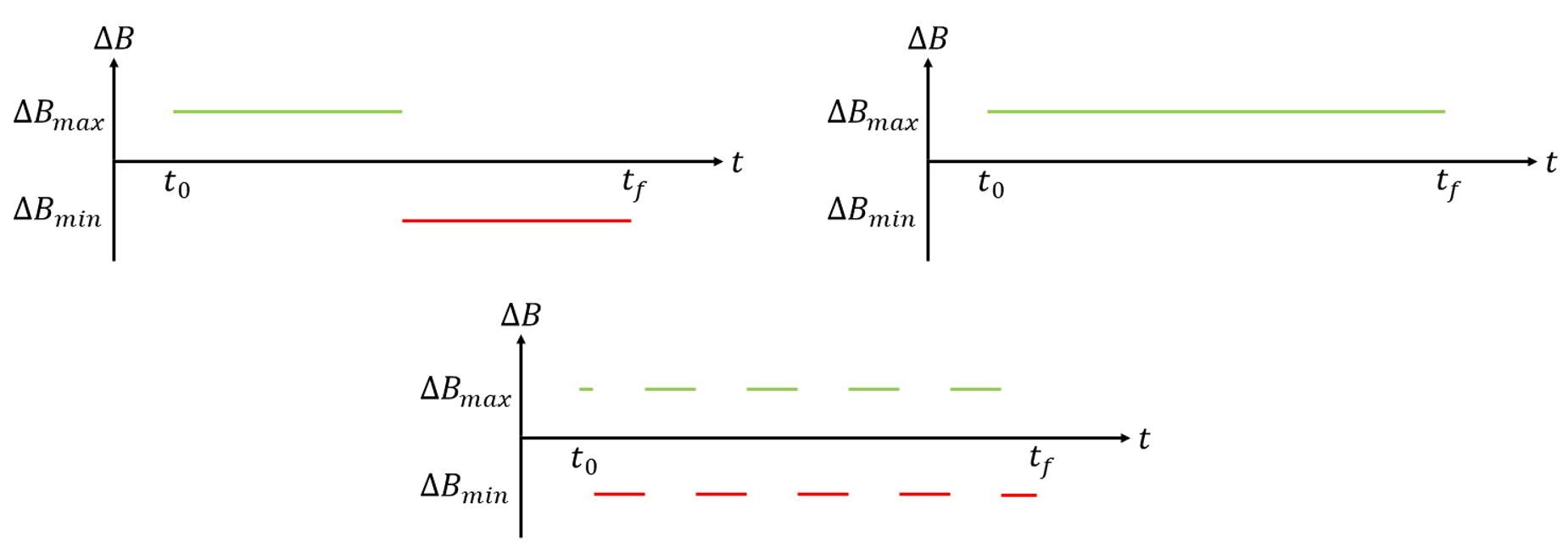}
    \caption{Optimal differential drag profile variations for strictly $\delta a$ (top left), $\delta \lambda$ (top right), and $\delta e$ (bottom) dominant reconfigurations \cite{hunter_hybrid}.}
    \label{fig: optimal drag profiles}
\end{figure}

Reachable set theory for continuous control identifies times for both the beginning and end of an optimal maneuver, effectively finding the optimal maneuver duration. In the $\Delta \delta \boldsymbol{a}$ plane, the $\delta a$ dominant maneuver spans the length of the control window in the direction of $\Delta \delta a$ (see $S^*$ contours matching $\eta_{\delta a}$ in Figure \ref{fig: drag reachable sets}, left). The $\delta \lambda$ dominant maneuvers span the first and second half of the control window in opposite directions (see $S^*$ contours matching $\eta_{\delta \lambda}$ in Figure \ref{fig: drag reachable sets}, left). Similar to propulsive control, the drag reachable set forms a consistently circular contour in the $\Delta \delta \boldsymbol{e}$ plane (see $S^*$ contours matching $\eta_{\delta e}$ in Figure \ref{fig: drag reachable sets}, right). $\delta e$ dominant maneuvers act maximally in the direction of $\Delta \delta \boldsymbol{e}$, such that they span the halfspace defined by $\eta$, and the optimal locations occur $\pm 90^\circ$ from the phase of $\Delta \delta \boldsymbol{e}$, at times given by the set
\begin{equation} \label{eq: T opt de drag}
    T_{\text{drag},\delta e} = \left( \phi_{\Delta \delta \boldsymbol{e}} - \pi/2 - \Dot{\omega}_c t_f - u_{c,0} + k_i\pi \right) / (n_c-\Dot{\omega}_c)
\end{equation}
\noindent Assuming $\Delta \delta \boldsymbol{\alpha}$ is not reachable with drag control alone, the optimal policy is always bang-bang control, as drag has no associated cost and is limited in magnitude by the bounds on $\Delta B^\rho$. Therefore, the drag profiles corresponding to each single dominance case are illustrated in Figure \ref{fig: optimal drag profiles}.

\section{HYBRID ARCHITECTURE}

The hybrid optimal solution architecture solves the full reconfiguration in three closed-form sequential steps: differential drag, OOP propulsive, and IP propulsive. The output from one section acts as an input into the next, such that $\Delta \delta \boldsymbol{\alpha}$ is modified by the expected effect of each maneuver plan and IP/OOP coupling is minimized. Restricting the hybrid control problem to a set of optimal maneuver times enables the construction of optimal hybrid maneuver plans in closed-form. For differential drag, the solver must find the optimal duration of each single dominance case optimal drag profile. In the interest of computational efficiency, the propulsive solver uses the smallest subset of optimal maneuver times that can achieve both optimality and controllability for IP and OOP states.

\subsection{Closed-form Differential Drag Solutions}

The goal of differential drag control is to find the attitude profile that minimizes the IP reconfiguration cost. The inclusion of zero-cost maneuvers into the hybrid control problem increases the likelihood that the reconfiguration will have multiple dominance cases with the same $\delta v_{\text{min}}$ (see $S^*$ contours labelled $\delta a/\lambda$ in Figure \ref{fig: drag reachable sets}, left). Therefore, the optimal drag maneuver plan varies more than the options shown in Figure \ref{fig: optimal drag profiles}. This complication is addressed by an optimal composite profile that combines portions of each single dominance case profile to accomplish the maximum IP cost reduction \cite{hunter_hybrid}. With a sufficiently long control window, this can result in the entire reconfiguration being reachable with drag control alone.

A key design choice in the implementation of the drag dynamics in Eq. \ref{eq: hybrid control problem} is the coarseness of the discretization of the control window. The LTV formulation approximates the continuous solution as $k$ increases to infinity. However, this will eventually require a prohibitively large amount of on-board memory, especially in the extended control windows used in hybrid control. To balance accuracy and computational efficiency, consider an approximation of the atmospheric density over the control window, given as
\begin{equation}
    \rho(t) =  \Bar{\rho} + at + b\sin\left( 2t/n_c \right)
\end{equation}
\noindent where $\Bar{\rho}$ is the time-averaged density over the control window, $at$ represents small secular changes over time, $b\sin\left( 2t/n_c \right)$ represents periodic changes due to the diurnal bulge near Earth's equator, and $a$ and $b$ are arbitrary coefficients. In the continuous case, the periodic effect disappears when integrated over an arbitrary half-orbit $[t,t+\pi/n_c]$. For a discretized control window, periodic effects can be removed prior to integration by sampling maneuvers spaced roughly $\pi/2n_c$ apart, such that the periodic terms between maneuvers cancel over the arbitrary half-orbit. Therefore, an appropriate coarse discretization is the set 
\begin{equation}
    T_{\text{drag}} = \left \{t_0, T_{\text{drag},\delta e}, T_{\text{drag},\delta e} + \frac{\pi/2}{n_c-\Dot{\omega}_c}, \frac{t_f-t_0}{2} \right \}
\end{equation}
\noindent This discretization allows the implementation of provably optimal control without significant loss in final state accuracy.

\subsection{Closed-form Propulsive Solutions}

The closed-form propulsive framework used by this paper expands on a previously developed approach for orbits of arbitrary eccentricity \cite{chernick}. First, calculate $\delta v_{\text{min}}$ for all dominance cases through Eq. \ref{eq: delta-v min} to determine the reconfiguration dominance case. Second, the relative effect $\Delta \delta \boldsymbol{\alpha}_i$ of each optimal maneuver time $t_i$ is found by calculating 
\begin{equation} \label{eq: manuever effect}
    \Delta \delta \boldsymbol{\alpha}_i = \boldsymbol{\Phi}^{J_2}(t_i)\boldsymbol{B}(t_i)\boldsymbol{u}^{\text{opt}}_i
\end{equation}
\noindent where $\boldsymbol{u}^{\text{opt}}_i$ is an optimal propulsive maneuver of magnitude $\delta v_{\text{min}}$. This effectively finds the set $S$ in Eq. \ref{eq: RST single maneuver}. Finally, all optimal maneuver plans are searched to find the combination with lowest cost. Linear systems of equations can be constructed with an optimality constraint on the non-dominant ROE, given as
\begin{align} \label{eq: optimal system}
    \begin{bmatrix}
        \Vec{\boldsymbol{1}} & \Vec{\Delta \delta \boldsymbol{\alpha}_i}
    \end{bmatrix}^{\text{T}}
    \Vec{\boldsymbol{c}} &= 
    \begin{bmatrix}
        1 & \Delta \delta \boldsymbol{\alpha}_i
    \end{bmatrix}^{\text{T}}
\end{align}
\noindent where $\Vec{\boldsymbol{c}} \in R^{m \times 1}$ is a column vector of magnitude coefficients for $m$ maneuvers, $\Vec{\boldsymbol{1}} \in R^{m \times 1}$ is a column vector of ones, and $\Vec{\Delta \delta \boldsymbol{\alpha}_i} \in R^{m \times m-1}$ concatenates values for the non-dominant ROE from Eq. \ref{eq: manuever effect}. Multiplying each coefficient in $\Vec{\boldsymbol{c}}$ by $\delta v_{\text{min}}$ provides the corresponding maneuver magnitudes.  For any combination of times that cannot produce a valid solution ($\sum_i c_i \neq 1$), a viable sub-optimal maneuver plan can be found by replacing $\Vec{\boldsymbol{1}}$ with the dominant ROE values from Eq. \ref{eq: manuever effect}.

Comparing the radial and tangential maneuver effects in Eq. \ref{eq: B matrix}, all optimal IP maneuvers are tangential-only. Consider restricting potential IP maneuver times to $T_{\text{opt},\delta e}$ in Eq. \ref{eq: T opt de prop}. Maneuvers at these times act in the same phase as $\Delta \delta \boldsymbol{e}$, simplifying the 2D $\Delta \delta \boldsymbol{e}$ problem to a 1D $|\Delta \delta \boldsymbol{e}|$ problem. The reduced IP problem is then fully controllable with 3 tangential maneuvers. The $\delta e$ dominance case is solved with non-dominant ROE $\Delta \delta a$ and $\Delta \delta \lambda$, and the $\delta a$ dominance case is solved with non-dominant ROE $|\Delta \delta \boldsymbol{e}|$ and $\Delta \delta \lambda$. The OOP problem can be solved with one or two normal maneuvers at maneuver times $T_{\text{opt},\delta i}$ in Eq. \ref{eq: T opt di}, with non-dominant ROE $\Delta \delta i_x$ as a placeholder for two-maneuver plans.

\subsection{Propulsive Accommodations for Hybrid Control}

The three sequential steps of the solution architecture do not mitigate the coupling effects of IP maneuvers on the OOP state, introducing the risk of significant OOP final state error. Given that all optimal IP maneuvers are tangential, the IP coupling effect can be found a-priori given the known resulting $\Delta \delta \lambda$ from these yet-unsolved tangential maneuvers $\boldsymbol{u}_{\text{T},k}$, given as
\begin{equation} \label{eq: OOP coupled dvmin}
    \Delta \delta i_{y,\text{OOP}} = \Delta \delta i_y - \sum_{k=1}^{p} \boldsymbol{\Gamma}(t_k) \boldsymbol{u}_{\text{T},k} = \Delta \delta i_y - \frac{\Phi_{61,t_0}^{J_2}}{\Phi_{21,t_0}^{J_2}} \Delta \delta \lambda
\end{equation}
\noindent where $\Phi_{21,t_0}^{J_2}$ and $\Phi_{61,t_0}^{J_2}$ correspond to $\boldsymbol{\Phi}^{J_2}(t = t_0)$.

\begin{figure}[h]
    \centering
    \includegraphics[scale=0.3]{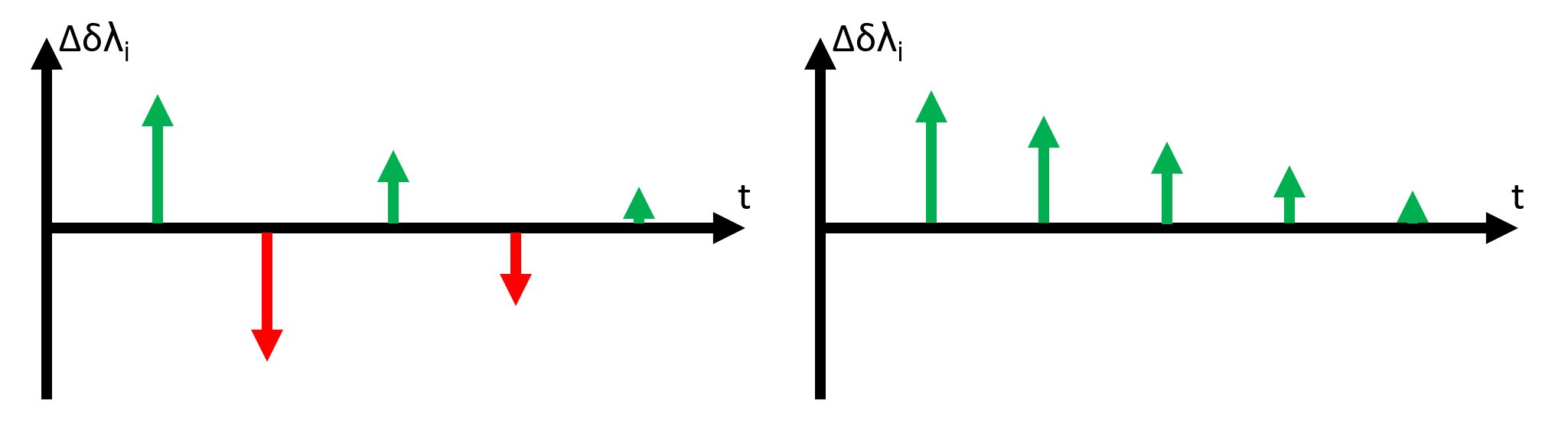}
    \caption{Relative $\Delta \delta \lambda$ effects for optimal maneuvers at times in $T_{\text{opt},\delta e}$ for the $\delta e$ dominance case (left) and the $\delta a$ dominance case (right).}
    \label{fig: dl effect}
\end{figure}

Searching through all possible IP maneuver plans requires $n$-choose-3 iterations, where $n$ is roughly double the number of orbits in the control window. Given the ${>}100$ orbit control windows used to maximize  differential drag, this is not a computationally viable approach. The computational load can be greatly reduced by identifying maneuvers that are guaranteed to be in a lowest cost maneuver plan. Considering the IP $\Delta \delta \boldsymbol{\alpha}_i$ from Eq. \ref{eq: manuever effect}, all maneuvers have the same $\Delta \delta a_i$ and $|\Delta \delta \boldsymbol{e}|_i$, but $\Delta \delta \lambda_i$ is time dependent and decreases in magnitude over the control window. Figure \ref{fig: dl effect} visualizes this general trend for both the $\delta e$ and $\delta a$ dominance cases to identify the most efficient maneuvers that span the $\Delta \delta \lambda$ state space. 

For the $\delta e$ dominance case, the first two times in $T_{\text{opt},\delta e}$ have the largest positive and negative $\Delta \delta \lambda_i$. For the $\delta a$ dominance case, the first and last time in $T_{\text{opt},\delta e}$ have the largest and smallest $\Delta \delta \lambda_i$. By restricting 3-maneuver plans to these maneuver times, searching over all plans only requires iterating through  $T_{\text{opt},\delta e}$ to find the third maneuver time. The algorithm now runs in linear time, approximately $2n$ iterations.

$\delta \lambda$ dominant reconfigurations do not fit the closed-form architecture presented thus far, as $\delta \lambda$ optimal maneuver times rarely align with the phase of $\Delta \delta \boldsymbol{e}$. A near-optimal solution can be generated by adding $\delta \lambda$ optimal maneuvers until the remaining control problem is $\delta e$ dominant. Consider a tangential maneuver at the beginning of the control window. The coefficient $c_0$ that causes the $\delta \lambda$ dominance case to equal the other IP dominance cases is given as
\begin{align}
    \delta v_{\text{min},\delta \lambda}(\Delta \delta \boldsymbol{\alpha} - c_0 \Delta \delta \boldsymbol{\alpha}_0) = \delta v_{\text{min, IP}}(\Delta \delta \boldsymbol{\alpha} - c_0 \Delta \delta \boldsymbol{\alpha}_0) \label{eq: dvmin dl = dvmin t0}
\end{align}
\noindent where $\delta v_{\text{min, IP}} = \delta v_{\text{min},\delta a}$ or $\delta v_{\text{min},\delta e}$, $\delta v_{\text{min}}(\Delta \delta \boldsymbol{\alpha})$ is the $\delta v_{\text{min}}$ for each dominance case as an analytical function of $\Delta \delta \boldsymbol{\alpha}$, and $\Delta \delta \boldsymbol{\alpha}_0$ is the pseudostate achieved by a unitary positive tangential maneuver at $t_0$. Using a positive unitary maneuver makes $c_0$ equivalently the size of the maneuver in units of m/s. If the $c_0$ for $\delta v_{\text{min},\delta \lambda} = \delta v_{\text{min},\delta a}$ has a smaller magnitude than that of $\delta v_{\text{min},\delta \lambda} = \delta v_{\text{min},\delta e}$, a $\delta a$-$\delta \lambda$ mutually optimal maneuver is added at the end of the control window. Eq. \ref{eq: dvmin dl = dvmin t0} is used again for $\delta v_{\text{min, IP}} = \delta v_{\text{min},\delta e}$, replacing $\Delta \delta \boldsymbol{\alpha}_0$ and $c_0$ with corresponding values at the end of the control window $\Delta \delta \boldsymbol{\alpha}_f$ and $c_f$. The previously presented closed-form architecture can then be run on the remaining $\delta e$ dominant control problem or skipped entirely if no $\delta \boldsymbol{e}$ control is desired.

\section{VALIDATION}

The effectiveness and validity of the proposed hybrid architecture is demonstrated through an example reconfiguration proposed for the upcoming SWARM-EX mission \cite{swarm-ex}. The SWARM-EX mission will use 3 identical 3U CubeSats, each with its own 3-axis attitude determination and control system (ADCS) and single thruster cold gas propulsion unit. The science mission will study the equatorial ionization anomaly (EIA) and equatorial thermospheric anomaly (ETA) in the ionized region of Earth's upper atmosphere. The mean orbital elements of the chief spacecraft, chosen arbitrarily among the 3 identical spacecraft, for this objective are given as
\begin{equation}
    \boldsymbol{\alpha}_c = 
    \begin{bmatrix}
        a & e & i & \Omega & \omega & M
    \end{bmatrix} = 
    \begin{bmatrix}
        6798\text{km} & 0.001 & 51^\circ & 0^\circ & 0^\circ & 90^\circ
    \end{bmatrix} 
\end{equation}
\noindent The SWARM-EX spacecraft have an estimated weight of 6kg, a modelled constant coefficient of drag of 1.5, and a minimum and maximum cross-sectional area dictated by its solar array and bounded between $0.01\text{m}^2$ and $0.09\text{m}^2$. Studying the full profile of the EIA and ETA requires about $1300\text{km}$ of inter-satellite separation at a tolerance of $\pm 100\text{km}$, while periodic instrumentation cross-calibration between satellites requires about $10\text{km}$ of separation while accounting for collision avoidance. These contrasting specifications define this "safe approach" validation reconfiguration for a 100 orbit control window, where the initial and final desired mean ROE state and resulting pseudostate of a deputy under J2 perturbations, dimensionally scaled by $a_c$, are defined in row vectors as

\begin{equation}
    \begin{bmatrix}
         a_c \delta \boldsymbol{\alpha}_0 \\ a_c \delta \boldsymbol{\alpha}_f \\
    \end{bmatrix} = 
    \begin{bmatrix}
        -0.02 & 100 & 0.5 & 0.9 & -0.05 & 0.95 \\
        0 & 10 & 0 & 1 & 0 & 1 \\
    \end{bmatrix} \text{ km}
\end{equation}
\begin{equation}
    a_c \Delta \delta \boldsymbol{\alpha} = 
    \begin{bmatrix}
        20 & -109,015 & -69 & -27 & 50 & 109
    \end{bmatrix} \text{ m}
\end{equation}

Three control algorithms will be used alongside the proposed hybrid architecture to evaluate its optimality, computational efficiency, and final state accuracy. The optimal hybrid control problem in Eq. \ref{eq: hybrid control problem} is convex and can be solved numerically with the interior point solvers in the CVX software package \cite{boyd}. CVX will verify that the differential drag solution optimally reduces the IP delta-v cost. Koenig and D'Amico created a reachable set-based impulsive control algorithm, referred to as the KD solver, to rapidly compute the lower bound for optimal control and numerically identify optimal maneuver times \cite{koenig_impulse}. The KD solver will quantify the error in the derived analytical $\delta v_{\text{min}}$ and define the expected final state error for a maneuver plan that uses J2-perturbed dynamics as the reference truth, neglecting other perturbations to Keplerian orbit. Chernick and D'Amico extended this approach to ROE-space to derive optimal closed-form solutions for most $\delta e$ and $\delta a$ dominant reconfigurations in orbits of arbitrary eccentricity, referred to as the CD solver \cite{chernick}. As mentioned previously, this solver forms the basis for the propulsive control algorithm used by this work, but does not included the accommodations for $\delta \lambda$ dominant maneuvering or extended control windows. The CD solver will demonstrate the improvements to provably optimal closed-form propulsive control presented in this work.

\subsection{Optimal Maneuver Planning}
Optimality is evaluated for each section of hybrid control (IP, OOP, drag) on the initial pseudostate, to separate it from accommodations for IP/OOP coupling. Considering only propulsive control, the proposed algorithm determines that the IP reconfiguration is $\delta \lambda$ dominant, with a $\delta v_{\text{min}}$ of 0.11894m/s. The KD algorithm finds an IP $\delta v_{\text{min}}$ of 0.11907m/s, a 0.13mm/s difference that is within the 1mm/s stopping tolerance set by the user. Therefore, the proposed algorithm is correctly identifying the optimal cost. The proposed algorithm results in an IP maneuver plan of cost 0.11938m/s, only 0.26\% sub-optimal. The CD algorithm also identifies the same dominance case and minimum cost as the proposed algorithm but produces an IP maneuver plan of cost 0.11989m/s (0.69\% sub-optimal), higher than the proposed algorithm.

For the OOP reconfiguration, the proposed algorithm finds the pseudostate in $\delta i$ dominance region 1 with a corresponding $\delta v_{\text{min}}$ of 0.10782m/s. The KD algorithm finds an OOP $\delta v_{\text{min}}$ of 0.10785/s, a 0.03mm/s difference that is within the stopping tolerance set by the user. Therefore, the proposed algorithm is correctly identifying the optimal cost. In contrast, the CD algorithm finds an OOP $\delta v_{\text{min}}$ of 0.13482/s, a 26.97mm/s error. Both the proposed algorithm and the CD algorithm produce an OOP maneuver plan with cost equal to its respective $\delta v_{\text{min}}$, such that the CD algorithm results in 25.01\% sub-optimality while the proposed algorithm outputs an optimal plan.

Now considering only differential drag control, the proposed algorithm finds an optimal drag profile corresponding to the strict $\delta \lambda$ dominance case and a remaining IP $\delta v_{\text{min}}$ of 0.0463m/s after applying the profile. This saves 0.07264m/s of delta-v cost (61.07\%) over propulsive-only methods, over half of the original cost. The numerical CVX interior point approach finds the same profile with a remaining IP $\delta v_{\text{min}}$ of 0.0465m/s, a minimal 0.02mm/s difference. It is worth noting here that the generic interior point solvers must have a coarse time discretization of propulsive and drag maneuver inputs to reliably converge to an optimal solution. This further demonstrates the advantages of implementing a numerically-tractable closed-form algorithm that can be easily tested for reliable control solution creation.


\subsection{Performance}

The performance and viability of CubeSat implementation is illustrated comparatively between different control approaches by evaluating computational efficiency and final state error of full maneuver plans. Computational efficiency is quantified by the run time required to produce a maneuver plan on equivalent hardware, in this case a 3.20 GHz processor and 16GB RAM. The run times required to solve the example reconfiguration for all considered control approaches are displayed in Table \ref{tab: alg run times}.

\begin{table}[h]
    \centering
    \caption{Algorithm run times and final state errors of various control approaches for the validation reconfiguration.}
    \label{tab: alg run times}
    \begin{tabular}{l|l|l}
        \hline
        Control Approach & Algorithm Run Time & $a_c\delta\boldsymbol{\alpha}_f$ Error \\
        \hline
        \hline
        Numerical Propulsive (KD) & 7.286s & 
        $\begin{bmatrix}
            -2 & 1672 & 4 & \hspace{4mm} 1 & \hspace{2mm} 1 & \hspace{2mm} -4
        \end{bmatrix}$ m \\
        \hline
        Closed-form Propulsive (CD) & 631.757s & 
        $\begin{bmatrix}
            -2 & 1643 & 4 & \hspace{1mm} -2 & -1 & 200
        \end{bmatrix}$ m \\
        \hline
        Numerical Hybrid (CVX)* & 6.887s & 
        $\begin{bmatrix}
            -2 & 1770 & \hspace{-1mm} 15 & \hspace{-1mm} -11 & \hspace{1mm} 2 & \hspace{2mm} -4
        \end{bmatrix}$ m \\
        \hline
        Closed-form Hybrid & 0.0012s & 
        $\begin{bmatrix}
            -2 & 1703 & 4 & \hspace{1mm} -2 & \hspace{2mm} 3 & \hspace{3mm} 25
        \end{bmatrix}$ m \\
        \hline
        \multicolumn{3}{c}{*Requires coarse time discretization to converge to solution}
    \end{tabular}
\end{table}

In general, closed-form control solvers are preferable for DSS, as the spacecraft often do not have the on-board computing power to store and process numerical optimization algorithms. However for the validation reconfiguration used here, the CD algorithm has the longest run time even though it is closed-form. This is because it searches through all possible 3-maneuver IP plans and therefore scales in polynomial time. The advances made by this work to identify efficient maneuver times a-priori allow the same approach to be solved almost instantaneously. It is worth noting here that the CVX solver having the same run time as the KD algorithm is not actually a product of superior efficiency but rather a sign of poor numerical stability. As mentioned before, a coarse time discretization is required for CVX to converge to the optimal solution of Eq. \ref{eq: hybrid control problem}, and this issue becomes even more apparent when considering the full 6D ROE reconfiguration. The KD algorithm solves the optimal impulsive control problem at the same time resolution regardless of control window length, at a much higher resolution than CVX can handle, and converges reliably to an optimal solution in a similar run time. This demonstrates the need for a specialized algorithm, numerical or closed-form, to be able to solve the hybrid control problem reliably and at high precision.

\begin{figure} [h]
    \centering
    \includegraphics[scale=0.5]{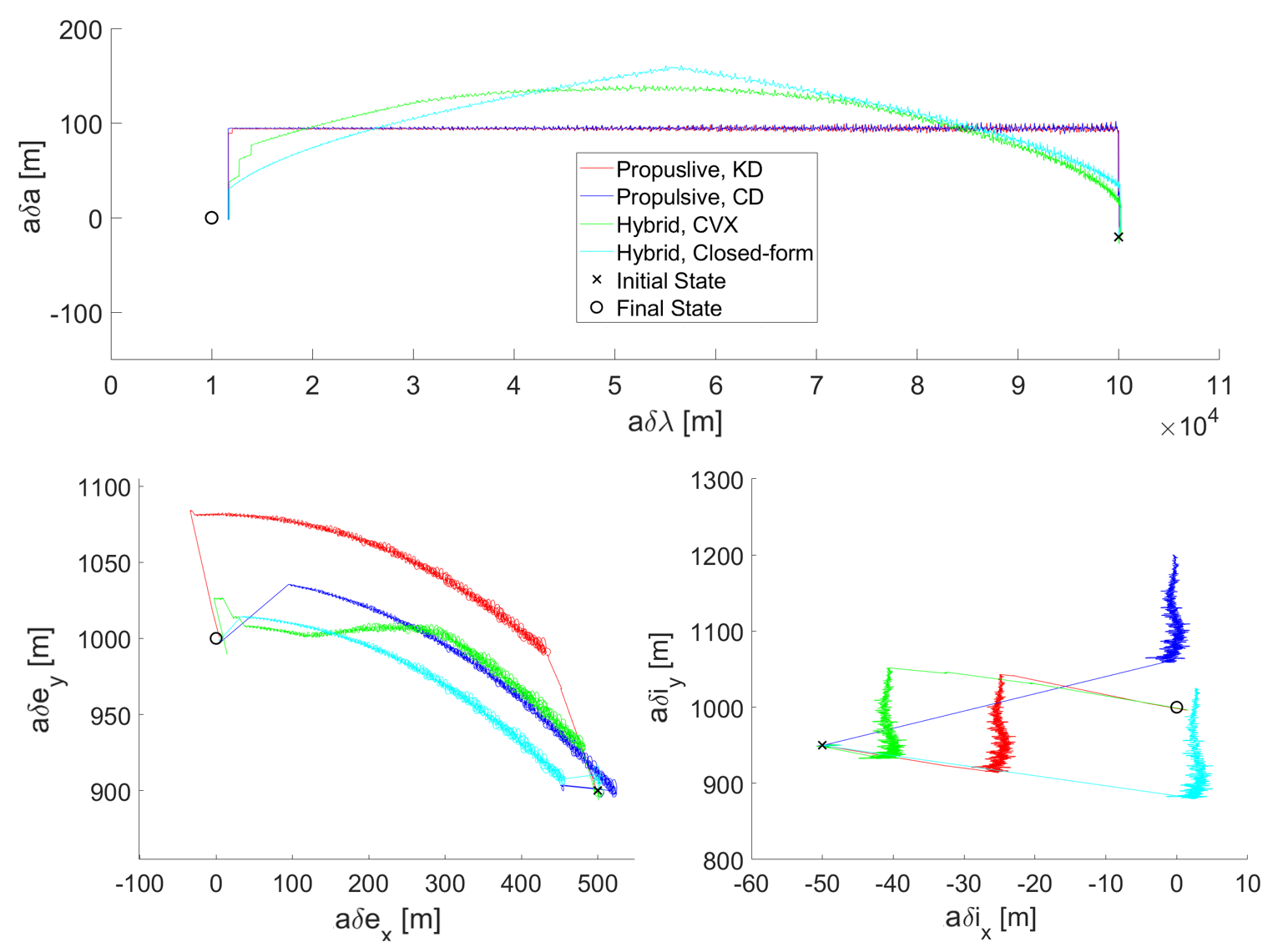}
    \caption{Mean ROE trajectories for the validation reconfiguration under a full-force simulation in the $\Delta \delta \boldsymbol{a}$ (top), $\Delta \delta \boldsymbol{e}$ (bottom left), and $\Delta \delta \boldsymbol{i}$ planes (bottom right).}
    \label{fig: control traj}
\end{figure}

All control solutions are propagated by a 4th-order Runge-Kutta integration of the GVE over a 10s time step, incorporating perturbations from Earth's gravity (30x30 gravity model), atmospheric drag using the NRLMSISE-00 density model, solar radiation pressure, third body effects, and relativistic corrections \cite{nrlmsise00}. Attitude constraints, such as Sun-pointing to remain power positive and Earth-pointing for atmospheric measurements, are neglected for simplicity as they do not affect the fundamental hybrid control problem and will only cause the minimum and maximum bounds on the differential attitude to become time dependent. Each spacecraft is assumed to have perfect knowledge of the atmospheric density on-board, to separate the validity of the hybrid approach from the problem of estimating and rejecting errors in the atmospheric density modelling. The ROE trajectory of the deputy spacecraft under each maneuver plan is displayed in Figure \ref{fig: control traj}.

The effect of the differential drag profile on hybrid control and corresponding delta-v savings is illustrated in the $\delta \boldsymbol{a}$ plane ($\delta a$, $\delta \lambda$). The parabolic trajectory caused by drag accomplishes a similar change in $\delta \lambda$ as propulsive-only approaches with lower magnitude tangential maneuvers at the beginning and end of the control window. Unlike the CD algorithm, the proposed closed-form hybrid approach is able to roughly match the final state error of the two numerical algorithms. The improvements made in this work on OOP maneuver plan accuracy are apparent in the $\delta \boldsymbol{i}$ plane, where the proposed algorithm has a much lower final state error than the CD approach. It also accomplishes this accuracy with fewer propulsive maneuvers than either numerical approach, demonstrated by the frequent instantaneous trajectory jumps for the numerical algorithms in the $\delta \boldsymbol{e}$ plane. The authors made a similar observation in differential drag maneuver frequency, in particular for $\delta \lambda$ dominant reconfigurations, implying that the overall proposed closed-form hybrid approach is more realizable under the true rotational dynamics of each DSS spacecraft \cite{hunter_hybrid}.

\section{CONCLUSIONS}

To take advantage of the unique capabilities of Distributed Space Systems (DSS), novel Guidance, Navigation, and Control (GNC) algorithms are required that enable spacecraft to conduct relative maneuvering at low propulsive cost and minimal computational load. Certain environmental forces, namely atmospheric drag in Low Earth Orbit (LEO), can be leveraged on-board to create differential forces with negligible propulsion. This paper proposes a full hybrid propulsive-differential drag architecture that offers significant delta-v savings over long, time-constrained control windows without sacrificing controllability. This architecture is realized through a Relative Orbital Element (ROE) state, simplifying nonlinear dynamics to a linear time variant (LTV) model, and reachable set theory, restricting the optimal control problem to a set of optimal maneuver times. The corresponding drag maneuver durations and propulsive magnitudes can be found in closed-form, producing provably optimal solutions for decoupled in-plane, out-of-plane ROE control. This deterministic methodology is shown to have comparable optimality and final state accuracy to numerical approaches while using minimal computing power. Errors in the dynamics model could potentially be mitigated through a Model Predictive Control (MPC) setup that re-solves the proposed algorithm at regular intervals, mimicking the iterative convergence properties found in the literature while maintaining provable optimality. Accommodations should be made to ensure that error rejection does not prohibitively increase delta-v cost, either by extending the control window such that the reconfiguration is reachable with drag control alone or by disabling propulsive control manually. The novel closed-form hybrid architecture can be used on a wide variety of DSS missions in LEO to increase mission lifetime. Furthermore, the general reachable set theory approach can be applied to any theoretical impulsive-continuous hybrid control problem with LTV dynamics to guide the creation of provably optimal closed-form solutions.

\section{Acknowledgment}
This work was supported by the National Science Foundation (Award Number 1936512). This project resulted from the Ideas Lab: Cross-cutting Initiative in CubeSat Innovations, an interdisciplinary program supported by Geosciences, Engineering, and Computer and Information Science and Engineering Directorates.

\bibliographystyle{AAS_publication}   
\bibliography{references.bib}   

\end{document}